\documentclass[a4paper,12pt]{amsproc}
\usepackage{amsmath,amsfonts,amssymb,amsthm,anysize}
\usepackage{enumerate}
\usepackage{epsfig}
\usepackage{supertabular}
\usepackage{graphicx}
\newcommand\N{{\mathbb N}}
\newcommand\Z{{\mathbb Z}}

\newcommand\F{{\mathbb F}}

\renewcommand\mod{{\mathrm{mod\, \, }}}

\theoremstyle{plain}
\newtheorem{theorem}{Theorem}[section]

\newtheorem{problem}[theorem]{Problem}
\newtheorem{lemma}[theorem]{Lemma}

\newtheorem{proposition}[theorem]{Proposition}
\newtheorem{example}[theorem]{Example}
\newtheorem{definition}[theorem]{Definition}

\numberwithin{equation}{section}

\theoremstyle{remark}
\newtheorem{remark}[theorem]{Remark}

\def\tra{{\mathrm{T}}}

\usepackage{url, hyperref}
\hypersetup{citecolor=blue, linkcolor=blue, colorlinks=true}

\renewcommand\ge{\geqslant}
\allowdisplaybreaks[2]
\begin{document}

\title[Difference families with four or eight blocks and Hadamard matrices]{New constructions of Hadamard matrices}

\author{Ka Hin Leung}
\address{ %
Department of Mathematics\\
National University of Singapore, 
Kent Ridge, Singapore 119260, Republic of Singapore}
\email{matlkh@nus.edu.sg}
\thanks{Ka Hin Leung was supported by Ministry of Education grant R-146-000-276-114}

\thanks{}

\author{Koji Momihara}
\address{ %
Division of Natural Science,\\
Faculty of Advanced Science and Technology,\\
Kumamoto University\\
2-40-1 Kurokami, Kumamoto 860-8555, Japan}
\email{momihara@educ.kumamoto-u.ac.jp}
\thanks{Koji Momihara was supported by 
JSPS under Grant-in-Aid for Young Scientists (B) 17K14236 and Scientific Research (B) 15H03636.}

\begin{abstract}
In this paper, we obtain a number of new infinite families of Hadamard matrices. Our constructions are based on 
four new constructions of difference families with four or eight blocks. By applying the Wallis-Whiteman array or the Kharaghani array to the difference families constructed, we obtain new Hadamard matrices of order $4(uv+1)$ for 
$u=2$ and  $v\in \Phi_1\cup \Phi_2 \cup \Phi_3 \cup \Phi_4$; and for  
$u\in \{3,5\}$ and $v\in \Phi_1\cup \Phi_2 \cup \Phi_3$. 
Here, $\Phi_1=\{q^2:q\equiv 1\,(\mod{4})\mbox{ is a prime power}\}$, $\Phi_2=\{n^4\in \N:n\equiv 1\,(\mod{2})\} \cup \{9n^4\in \N:n\equiv 1\,(\mod{2})\}$, $\Phi_3=\{5\}$ and $\Phi_4=\{13,37\}$. Moreover, our construction also yields new Hadamard matrices of order $8(uv+1)$ for any  $u\in \Phi_1\cup \Phi_2$ and $v\in \Phi_1\cup \Phi_2 \cup \Phi_3$.

\end{abstract}



\maketitle

\section{Introduction}
Throughout this paper, we fix the 
following notation. Let $(G, +)$ be an additively written abelian group of order $v$ and let $G^\ast:=G\setminus \{0_G\}$. For any subset $D$ in $G$, we denote 
$D^{(-1)}:=\{-x:x \in D\}$, 
$\overline{D}:=G^\ast \setminus D$, and $D^c:=G \setminus D$. 
Furthermore, we also denote $\sum_{x\in D}x\in \Z[G]$ by $D$ when there is no confusion.  

\medskip

Let $B_i$, $i=1,2,\ldots,\ell$, be $k_i$-subsets of $G$ and ${\mathcal B}=\{B_i:i=1,2,\ldots,\ell\}$. 
Any subset $A$ of $G$ is said to be {\it symmetric} if $ A=A^{(-1)}$.  A family ${\mathcal B}$ is said to be {\it symmetric} if each $B_i$ is symmetric.  ${\mathcal B}$ is said to be a 
{\it difference family with parameters $(v;k_1,k_2,\ldots,k_\ell;\lambda)$ in $G$} if 
the list of differences ``$x-y, x,y\in B_i,x\not=y,i=1,2,\ldots,\ell$" represents every element of $G\setminus \{0_G\}$ exactly $\lambda$ times; or equivalently  
\[ \sum_{i=1}^{\ell} B_i B_i^{(-1)}= \lambda G+ \Big(\sum_{i=1}^{\ell} k_i-\lambda\Big) \cdot 0_G.\]

\begin{remark}
To avoid possible confusion, we write $\lambda \cdot 0_G$ as the element in the group ring such that the coefficient of $0_G$ is $\lambda$. This notation will be used throughout this paper. 
\end{remark}

Each subset $B_i$ is called 
a {\it block} of ${\mathcal B}$. If there is only one block in a difference family, the block is a {\it difference set}. 

\medskip

In this paper, our main idea is to construct a new difference family in $G\times G'$ from difference families in $G$. Thus, we need the difference families in $G$ to satisfies certain properties. We first define the following:

\begin{definition}
Suppose ${\mathcal B}=\{B_i:i=1,2,\ldots,\ell \}$ is  a  difference family in $G$.

(i) For $\ell=2,4,8$, ${\mathcal B}$ is  said to be of 
{\it type $H_\ell^\ast$} 
if  $\sum_{i=1}^\ell k_i-\ell (|G|+1)/4=\lambda$. 

(ii) For $\ell=4$, ${\mathcal B}$ is said to be of {\it type $H$} if $\sum_{i=1}^4 k_i-|G|=\lambda$. 
\end{definition}

Let us first explain how a difference family can be used to construct a Hadamard matrix. 
It is well known that if there is a difference family of type $H$ in $G$, then we have a Hadamard matrix of order $4|G|$ by plugging the circulant $(-1,1)$ matrices obtained from its blocks into 
the Goethals-Seidel array~\cite{GS67}. In literature, a  number of difference families of type $H$ have been extensively studied~\cite{LMS06,XL91,XXSW05,X98,Y92}.

For any difference family of type $H_\ell^\ast$ in $G$ for $\ell=2$ or $4$, we construct a Hadamard matrix of order $\ell(|G|+1)$ by plugging the circulant $(-1,1)$ matrices obtained from its blocks into 
the Szekeres array~\cite[Theorem~2.4]{WSW} or the Wallis-Whiteman array~\cite[Theorem~4.17]{WSW}, respectively. In the case where $\ell=8$, we may use the Kharaghani array under the assumption that the circulant $(-1,1)$ matrices $M_i$, $i=0,1,\ldots,7$, obtained from its blocks are amicable, i.e., $\sum_{i=0}^3
(M_{2i}M_{2i+1}^\tra-M_{2i+1}M_{2i}^\tra)=O$. (See, e.g., 
\cite{K00},\cite[Lemma~4.20]{S17},\cite[Page~12]{SY89}.)  Note that the amicability condition becomes trivial if the difference family is symmetric. That explains why we are interested in difference families of type $H_4^\ast$ and symmetric difference families of type $H_8^\ast$ in this paper. 

\medskip

As mentioned before, if a difference family consists of only one block, the block is actually a difference set. However,  the parameters of the single block must be of the form $(4n^2;2n^2-n;n^2-n)$. 
In the literature, a difference set with parameters $(4n^2;2n^2-n;n^2-n)$ is called a {\it Menon difference set}. For any Menon difference set  in a group of order $4n^2$, we consider  an incidence matrix of the design developed; by replacing all zero entries with $-1$, we obtain a Hadamard matrix of order $4n^2$. The existence of Menon difference sets in abelian groups of order $4n^2$ with $n$ odd have been extensively studied~\cite{C97,T84,X92,WX97,XC96}. It is known that for any  even $s\ge 0$ and $k\ge 1$, a Menon difference set  in a group of order $4\cdot 2^s k^4 $ or $4\cdot 9\cdot 2^s k^4$ exists. 
Moreover, Turyn~\cite{T84} has shown how to construct a Menon difference set  in a group of
order $4u^2n^2$ from Menon difference sets  in groups of  order $4u^2$ and $4n^2$. It is then natural to ask if we can extend those results to difference families. It turns out that we need some extra conditions imposed on the differences families concerned. We first list the conditions required. 

\medskip 

\begin{definition}
Let ${\mathcal D}=\{D_i:i=0,1,2,3\}$ be a difference family of type $H$ in an abelian group $G$.   We define conditions (c1)-(c5) as follows:  
\begin{itemize}
\item[(c1)] $D_0D_2^{(-1)}+D_2D_0^{(-1)}+D_1D_3^{(-1)}+D_3D_1^{(-1)}=
(\sum_{i=0}^3|D_i|-|G|)G$.  
\item[(c2)] For each $i=0,1,2,3$, $D_i$ is symmetric; and 
\[ D_0D_2^{(-1)}+D_2D_0^{(-1)}+D_1D_3^{(-1)}+D_3D_1^{(-1)}=
\Big(\sum_{i=0}^3|D_i|-|G|\Big)G.\]
\item[(c3)] For each $i=0,1,2,3$, $D_i$ is symmetric; and for $(h,i,j,k)=(0,1,2,3), (0,2,1,3),$ $(0,3,1,2)$,
\[ D_hD_i^{(-1)}+D_iD_h^{(-1)}+D_jD_k^{(-1)}+D_kD_j^{(-1)}=
\Big(\sum_{i=0}^3|D_i|-|G|\Big)G.\]
\item[(c4)] $2 | (\sum_{i=0}^3D_i)$ in $\Z[G]$, i.e. all coefficients in $\sum_{i=0}^3D_i$ are even. 
\item[(c5)] $D_0D_2^{(-1)}=D_2D_0^{(-1)}$ and $D_1D_3^{(-1)}=D_3D_1^{(-1)}$.  
\end{itemize}
\end{definition}
\medskip

It has been shown in \cite{XX99} that  a difference family of 
type $H$ satisfying (c4) yields a ``$T$-matrix'', which has played an important role in constructions of Hadamard matrices \cite{T74,XX99}.
It is also known that a difference family of type $H$ satisfying  (c1) yields a difference family with parameters $(2v;k_1,k_2;k_1+k_2-v)$ \cite{XL95}.  We may interpret that such a difference family with two blocks belongs to the class which is in between difference families of type $H$  and Menon difference sets. 
First, the number of blocks for such a family is $2$ whereas the number of blocks in a difference family of type $H$ and a Menon difference sets are $4$ and $1$ respectively. Second, a 
difference family with parameters $(2v;k_1,k_2;k_1+k_2-v)$ can be viewed as  subsets of a group of order $2v$ whereas a difference family of type $H$ can be viewed as subsets of a group of order $4v$ and a Menon difference set sits in a group of order $v$.

\medskip

 The following `product constructions' using difference families of type $H$ were first studied in \cite{XL95,XX94,XX99}.  

\begin{theorem}\label{thm:consttypeH}
Let ${\mathcal D}$ and ${\mathcal D'}$ be difference families of type $H$ in an abelian group of order $u$ and $n$ respectively. If any of the following conditions is satisfied, then there exists  a difference family of type $H$ in an abelian group of order $un$. 
\begin{itemize}
\item[(1)]~\cite{XL95} ${\mathcal D}$ satisfies  (c1) and  ${\mathcal D'}$ satisfies (c3).
 \item[(2)]~\cite{XL95}  ${\mathcal D}$ satisfies  (c1) and  ${\mathcal D'}$ satisfies (c2). 
\item[(3)]~\cite{XX99}    ${\mathcal D}$ satisfies  (c4) and  ${\mathcal D'}$ satisfies (c5).
\item[(4)]~\cite{XX99}   ${\mathcal D}$ satisfies  (c4) and  ${\mathcal D'}$ satisfies (c1). 
\end{itemize}
Moreover, if (1) is satisfied, then there exists a difference family of type $H$ that satisfies (c1). 
\end{theorem}
\begin{remark}\label{rem:intro}
Note that Theorem~\ref{thm:consttypeH} (2) is not stated explicitly in \cite{XL95}. But it can be proved easily by using a similar argument as in the proof of Theorem~2 in \cite{XL95}. 
\end{remark}

\medskip

Unfortunately, Theorem~\ref{thm:consttypeH} is not always applicable. A somewhat different  `product construction'  for ``Paley type partial difference sets'' has been considered by J.~Polhill~\cite{P10}. A {\it Paley type partial difference set} is a symmetric subset $D$ of size $(|G|-1)/2$  in  a finite group $G$  satisfying 
$D^2=\frac{|G|-1}{2} \cdot 0_G-D+\frac{|G|-1}{4}G^\ast$. 
To construct Paley type partial difference sets in $G\times G'$, it is not sufficient just to make use of one Paley type partial difference sets in $G$ and one in $G'$. Indeed, in Theorem 4.1 ~\cite{P10}, one needs a difference family of type $H$ and a Paley partial difference set in each group $G$ and $G'$.  
On the other hand, if $P$ is a Paley type partial difference set, then $\{P, \overline{P}\}$ is a symmetric difference family of type $H_2^\ast$. Therefore, his construction can be viewed as a  `product construction' that involves two difference families, where one is of type $H$ and the other is of type $H_2^\ast$. 
It is therefore natural to ask the following:

\begin{problem}
Is it possible to construct a difference family of type $H_4^\ast$ or $H_8^\ast$ in the group $G\times G'$ using two difference families of type $H$ and $H_4^\ast$ 
in  each $G$ and $G'$?
\end{problem} 

To address  the problem above, we need to work with difference families that satisfy certain conditions.  

\begin{definition} \label{def:E}
Let  ${\mathcal D}=\{D_i:i=0,1,2,3\}$ be a  difference family of  type $H$ in a group $G$. We define conditions (d1)-(d5) as follows:
\begin{itemize}
\item[(d1)] $0_G\not \in \bigcup_{i=0}^3D_i$. 
\item[(d2)] There are subsets $E_0,E_1$ of $G^\ast$ such that 
$D_0+D_1-D_2-D_3=E_0-\overline{E}_0$, $D_0+D_3-D_1-D_2=E_1-\overline{E}_1$, and the set  of $E_0,E_1,\overline{E}_0,\overline{E}_1$ forms a difference family of type $H_4^\ast$. 
\item[(d3)] $D_0D_2^{(-1)}+D_2D_0^{(-1)}+D_1D_3^{(-1)}+D_3D_1^{(-1)}=(\sum_{i=0}^3|D_i|-|G|+1)G$.
\item[(d4)] The coefficient of each element $x\in G^\ast$ in $\sum_{i=0}^3D_i$ is one or three, and the coefficient of $0_G$ in $\sum_{i=0}^3D_i$ is even but not two in exactly one of
$\sum_{i=0,1}D_i$, $\sum_{i=2,3}D_i$, $\sum_{i=0,3}D_i$ and 
$\sum_{i=1,2}D_i$. 
\item[(d5)] $|D_0|+|D_3|=|D_1|+|D_2|$. 
\end{itemize}
\end{definition}

Our first goal is to show the following:

\begin{theorem}\label{thm:consttypeH2}
Let ${\mathcal D}=\{D_i:i=0,1,2,3\}$  be a difference family of type $H$ in a group $G$. 
\begin{itemize}
\item[(1)] If ${\mathcal D}$  satisfies  (d3);  and there exists a difference family of type $H_4^\ast$ satisfies (c1) in $G$, then there exists a difference family of type $H_4^\ast$ in $G\times \Z_2$.  
\item[(2)] If ${\mathcal D}$  is symmetric and it satisfies  (d3) and there exists a difference family of type $H_2^\ast$ in $G$, then there exists a difference family of type $H_4^\ast$ in $G\times \Z_3$. 
\item[(3)]  If ${\mathcal D}$ is symmetric and it satisfies  (d2) and (d5),  then there exists a difference family of type $H_4^\ast$ in $G\times \Z_5$.
\item[(4)]  If ${\mathcal D}$ is symmetric and it satisfies  (d2);  and  there exists a symmetric difference family ${\mathcal S}$ of type $H$  satisfying  (d1) and (d2) in an abelian group $N$, then there exists a  symmetric difference family of type $H_8^\ast$  in an abelian group of order $G\times N$. 
\end{itemize}
\end{theorem}

Before we will prove Theorem~\ref{thm:consttypeH2} in Section~\ref{sec:const},  
we will first study some general properties of difference families of type $H$ in Section~\ref{sec:proper}. In  Section~\ref{sec:existsDF}, 
we construct explicitly difference families of type $H$ that can be used in Theorem~\ref{thm:consttypeH2}. In the final section, we prove our main results concerning the existence of Hadamard matrices. 

\begin{theorem}\label{thm:existHada}
Let $\Phi_1=\{q^2:q\equiv 1\,(\mod{4})\mbox{ is a prime power}\}$, $\Phi_2=\{n^4\in \N:n\equiv 1\,(\mod{2})\} \cup \{9n^4\in \N:n\equiv 1\,(\mod{2})\}$, $\Phi_3=\{5\}$ and $\Phi_4=\{13,37\}$. Then, the following hold: 
\begin{itemize}
\item[(1)]
There exists a Hadamard matrix of order $4(2v+1)$ for $v\in  \Phi_1 \cup \Phi_2 \cup \Phi_3 \cup \Phi_4$.
\item[(2)]
There exists a Hadamard matrix of order $4(3v+1)$ for $v\in  \Phi_1 \cup \Phi_2 \cup \Phi_3$. 
\item[(3)] There exists a Hadamard matrix of order $4(5v+1)$ for 
$v\in \Phi_1 \cup \Phi_2 \cup \Phi_3$. 
\item[(4)] There exists a Hadamard matrix of order $8(uv+1)$ for $u\in \Phi_1 \cup \Phi_2$ and $v\in \Phi_1 \cup \Phi_2 \cup \Phi_3$. 
\end{itemize}
\end{theorem}

It is generally difficult to construct 
Hadamard matrices of order $4p$ for an odd prime $p$.  
Note that all elements in $\Phi_1 \cup \Phi_2 \cup \Phi_3 \cup \Phi_4$ are odd. Thus, we may only apply
Theorem~\ref{thm:existHada}~(1) to generate Hadamard matrices of order $4p$ where $p$ is an odd prime.  
If $v\in \Phi_1$, then $3|(2v+1)$. Hence, we consider only those $v$ in $\Phi_2$, i.e. $v=n^4$ or $9\cdot n^4$ for an odd integer $n$. 
For odd  $n<300$, $2n^4+1$ is a prime for the following $n$, 
\[
1,
3,
21,
45,
63,
81,
105,
153,
177,
201,
219,
225,
249,
279,
297.
\]
Indeed, for odd $n<1000$, there are altogether $32$ such primes. 

On the other hand, for odd $n<1000$, there are altogether $74$ primes of the form 
$2\cdot 9 n^4+1$. For odd $n<300$, $2\cdot 9 n^4+1$ is a prime for the following $n$:
\begin{align*}
&1,
3,
5,
31,
45,
55,
57,
71,
79,
89,
107,
109,
119,
123,
137,
\\
&141,
159,
167,
173,
181,
197,
217,
255,
275,
285,
295.
\end{align*}

At this point, we are not able to prove if the sets $\{2n^4+1 : n \mbox{ is odd}\}$ and 
$\{2\cdot 9 n^4+1  : n \mbox{ is odd}\}$ contain infinitely many primes though.

\section{Symmetric difference families of type $H$ that satisfy (d2)}\label{sec:proper}

In this section, we aim to study  a symmetric difference family 
of type $H$ that satisfies (d2). Let us first give an example. 

\begin{example}\label{exam:v5}
Let $v=5$.  Define 
\[
D_0=\{0\},D_1=\{1,4\}, D_2=\{0\},D_3=\{2,3\}. 
\]
Then, the family ${\mathcal D}=\{D_i:i=0,1,2,3\}$ is a symmetric difference family of type $H$  in $\Z_5$ that satisfies (d2).  Note that ${\mathcal D}$ also satisfies (d5) but not (d1). 
\end{example}

In view of the above example, it is then interesting to explore how condition (d2) and other conditions are related. It is straightforward to check the following: 

\begin{proposition}\label{prop:comp}
Let ${\mathcal D}=\{D_i:i=0,1,2,3\}$ be a symmetric difference family of type $H$ satisfying (d2) in $G$ and that $0_G\in \bigcap_{i=0}^3 D_i$. Then ${\mathcal D}'=\{G\setminus D_i:i=0,1,2,3\}$ is a symmetric difference family of type $H$ that satisfies (d1) and (d2).   
\end{proposition}

From now on, we always assume ${\mathcal D}=\{D_i:i=0,1,2,3\}$ is a symmetric difference family of type $H$ satisfying (d2) in $G$; and we define $E_i$'s as defined in Definition~\ref{def:E}. To simplify our notation, we set
\begin{align*}
W&\,:= D_0D_2+D_1D_3;\\
X_0&\,:= D_0+D_1-D_2-D_3;\\
X_1&\,:= D_0+D_3-D_1-D_2;\\
T_0&\,:= (G+D_2-D_0)E_0+
(G-D_2+D_0)\overline{E}_0+(G+D_1-D_3)E_1+
(G-D_1+D_3)\overline{E}_1;\\
T_1&\,:=(G+D_3-D_1)E_0+
(G-D_3+D_1)\overline{E}_0+(G+D_2-D_0)E_1+
(G-D_2+D_0)\overline{E}_1;\\
T_2&\,:=(G+D_0-D_2)E_0+
(G-D_0+D_2)\overline{E}_0 +(G+D_3-D_1)E_1+
(G-D_3+D_1)\overline{E}_1;\\
T_3&\,:=(G+D_1-D_3)E_0+
(G-D_1+D_3)\overline{E}_0+(G+D_0-D_2)E_1+
(G-D_0+D_2)\overline{E}_1
\end{align*}
throughout this section. 

\begin{proposition}\label{prop:key1}
Suppose ${\mathcal D}$ is symmetric and of type $H$. Then ${\mathcal D}$ satisfies (d2) if and only if   ${\mathcal D}$ satisfies (d3) and (d4).  
\end{proposition}
\proof  Suppose ${\mathcal D}$ satisfies (d2). Then $X_0=E_0-\overline{E}_0$,  $X_1=E_1-\overline{E}_1$, and the set  of $E_0,E_1,\overline{E}_0,\overline{E}_1$ forms a difference family of type $H_4^\ast$. 
It is clear that ${\mathcal D}$ satisfies (d4). 
Observe that all $D_i$'s are symmetric. Hence, 
\begin{equation}\label{eq:w1}
 4W    = 2\sum_{i=0}^3 D_i^2 -X_0^2 -X_1^2 = 2\sum_{i=0}^3 D_i^2 - (E_0-\overline{E}_0)^2- (E_1-\overline{E}_1)^2 . \end{equation}
By assumption, ${\mathcal D}$ is of type $H$, therefore, 
\begin{equation} \label{eq:w2}
\sum_{i=0}^3 D_i^2= \Big(\sum_{i=0}^3|D_i|\Big) \cdot 0_G+
\Big(\sum_{i=0}^3|D_i|-|G|\Big)G^\ast. \end{equation}
On the other hand, 
\begin{align}
  (E_0-\overline{E}_0)^2+ (E_1-\overline{E}_1)^2 & = \sum_{i=0,1}(E_i^2+\overline{E}_i^2)-
 2 \sum_{i=0,1} E_i\overline{E}_i \nonumber\\
     &  = \sum_{i=0,1}(E_i^2+\overline{E}_i^2)+\Big(\sum_{i=0,1}(E_i+ \overline{E_i})^2
-\sum_{i=0,1}(E_i^2+\overline{E}_i^2)\Big)\nonumber\\
 & = 2(|G|-1)\cdot 0_G-2G^\ast. \label{eq:w3}
\end{align}
Combining \eqref{eq:w1},  \eqref{eq:w2} and \eqref{eq:w3}, it is easy to check that  (d3) is satisfied.

Conversely, we assume ${\mathcal D}$ satisfies (d3) and (d4).  By (d3), we have 
\begin{equation}\label{eq:si2}
\sum_{i=0,1}X_i^2
=2\sum_{i=0}^3 D_i^2-4W
=(2|G|-2) \cdot 0_G-2G^\ast. 
\end{equation}
 We consider the coefficients $a_x, b_x$ of any $x\in G$ in  $D_0+D_1-D_2-D_3$ and  $D_0+D_3-D_1-D_2$ respectively.  
 By condition (d4),  if $x=0_G$,  $a_x=b_x=0$.  If $x\neq 0_G$, then $a_x=\pm 1$ and $b_x=\pm 1$. 
  Thus, there are symmetric subsets $E_0,E_1$ of $G^\ast$ such that 
$X=E_0-\overline{E}_0$ and $X_1=E_1-\overline{E}_1$. Furthermore, we have 
\begin{equation}\label{eq:si3}
\sum_{i=0,1}X_i^2=2\sum_{i=0,1}(E_i^2+\overline{E}_i^2)-2\cdot 0_G-2(|G|-2)G. 
\end{equation}
Then, combining \eqref{eq:si2} and \eqref{eq:si3}, we have
\[ \sum_{i=0,1}(E_i^2+\overline{E}_i^2)=(2|G|-2)\cdot 0_G+(|G|-3)G^\ast.\]
Therefore, 
the set of $E_0,E_1,\overline{E}_0,\overline{E}_1$ forms a difference family of type $H_4^\ast$ in $G$. 
\qed

\begin{lemma}\label{lem:key2}
Let ${\mathcal D}$ be a symmetric difference family of type $H$  that satisfies (d2) in $G$. Then
\[ T_0=T_1=(|G|-1)\cdot 0_G+(2|G|-1)G^\ast  \mbox{ and } T_2=T_3=3(|G|-1)\cdot 0_G+(2|G|-3)G^\ast.\]
\end{lemma}
\proof  The proof for obtaining all $T_i$'s are similar. Here, we just work with $T_0$. 
By  (d2), we see that 
\[ D_0+D_1-D_2-D_3=E_0-\overline{E_0} \mbox{ and } D_0+D_3-D_1-D_2=E_1-\overline{E_1}.\]
Therefore, 
\begin{align}
T_0&\,=(|E_0|+|E_1|+|\overline{E_0}|+|\overline{E_1}|)G+(D_2-D_0)(E_0-\overline{E_0})+(D_1-D_3)(E_1-\overline{E_1})
\nonumber\\
&\,=2(|G|-1)G-(D_0-D_2)(D_0+D_1-D_2-D_3)-(D_3-D_1)(D_0+D_3-D_1-D_2) \nonumber\\
&\,=2(|G|-1)G+2W-\sum_{i=0}^3D_i^2=(|G|-1)\cdot 0_G+(2|G|-1)G^\ast.\label{eq:h1}
\end{align}
Note that the last equality in \eqref{eq:h1} follows from Proposition~\ref{prop:key1} and the assumption that ${\mathcal D}$ is of type $H$. 
\qed

\medskip

It turns out that  a  symmetric difference family of type $H$ that satisfies (d2) can be constructed using a family of symmetric subsets that satisfies certain conditions. 

\begin{definition}
We say that a family ${\mathcal A}=\{A_i:i=0,1,\ldots,7\}$ of symmetric subsets of $G$ is a {\bf building family} if it satisfies the following conditions: 
 \begin{itemize}
\item[(a1)] $A_i\cap A_j=\emptyset$ for $i\not =j$;
\item[(a2)] $\bigcup_{i=0}^7A_i=G\setminus \{0_G\}$; 
\item[(a3)] the set of $\bigcup_{i=0,1,6,7}A_i$, $\bigcup_{i=2,3,4,5}A_i$, 
$\bigcup_{i=0,3,5,6}A_i$, $\bigcup_{i=1,2,4,7}A_i$ forms a 
difference family of type $H_4^\ast$ in $G$; 
\item[(a4)] $\sum_{i=0}^7A_i^2-\sum_{i=0}^7A_iA_{i+4}
=(|G|-1) \cdot 0_G+\sum_{i=0}^3A_i-\sum_{i=4}^7A_i$, where the subscript of 
$A_{i+4}$ is reduced modulo $8$. 
\end{itemize}
\end{definition}

\begin{example}
Let $G=\Z_3\times \Z_3$. Define 
\[
A_0=\{(0,1),(0,2)\},A_1=\{(1,0),(2,0)\},A_2=\{(2,1),(1,2)\},
A_3=\{(1,1),(2,2)\}
\]
and $A_4=A_5=A_6=A_7=\emptyset$. 
Then, the family ${\mathcal A}=\{A_i:i=0,1,\ldots,7\}$ is a building family.  
\end{example}

Davis and Jedwab ~\cite{DJ} defined the notion of covering extended building sets in an abelian group and they are used to construct a difference set. Here, a building family in $G$ turns out to be 
a key component in our constrution of symmetric difference family of type $H$.

\begin{proposition}\label{prop:equiv}
Let ${\mathcal A}=\{A_i:i=0,1,\ldots,7\}$ and ${\mathcal D}=\{ D_0,D_1,D_2,D_3\}$ be families of symmetric subsets in $G$ such that  
\[ D_0= A_4 \cup A_1 \cup A_2 \cup A_3, \ D_1= A_5 \cup A_0 \cup A_2 \cup A_3,  \ D_2= A_6 \cup A_0 \cup A_1 \cup A_3, \ D_3=A_7 \cup A_0 \cup A_1 \cup A_2. \]
Then ${\mathcal A}$ is a building family  if and only if 
${\mathcal D}$ is a symmetric difference family of type $H$ 
that satisfies (d1) and (d2). 
\end{proposition}
\proof Let ${\mathcal A}$ and ${\mathcal D}$ be as defined above. First, we prove two claims. 

{\it Claim 1.} Suppose ${\mathcal A}$ satisfies (a2). Then ${\mathcal D}=\{D_i:i=0,1,2,3\}$ forms a symmetric difference family of type $H$ if and only if (a4) is satisfied. 

Observe that 
\begin{align}
\sum_{i=0}^3D_i^2=&\,\sum_{i=4}^7A_i^2+3\sum_{i=0}^3A_i^2
+2\sum_{i,j=0:i\not=j}^3 A_iA_j+2\sum_{i,j=0;i\not=j}^3A_iA_{j+4}\nonumber\\
=&\,\sum_{i=0}^7A_i^2 +2\Big(\sum_{i=0}^3A_i\Big)^2+2\sum_{i,j=0;i\not=j}^3A_iA_{j+4}. \label{Ka:h2}
\end{align}
On the other hand, 
\begin{align}
 &\,[\sum_{i=0}^7A_iA_{i+4}+(|G|-1)\cdot 0_G+\sum_{i=0}^3A_i-\sum_{i=4}^7 A_i]
+2\Big(\sum_{i=0}^3A_i\Big)^2+2\sum_{i,j=0;i\not=j}^3A_iA_{j+4}\label{Ka:h3} \\ 
=&\,(|G|-1)\cdot 0_G+\sum_{i=0}^3A_i-\sum_{i=4}^7 A_i
+2\Big(\sum_{i=0}^3A_i\Big)G^\ast \nonumber\\
= & \,(|G|-1)\cdot 0_G+ 2\Big(\sum_{i=0}^3|A_i|\Big)G - \sum_{i=0}^7 A_i \nonumber \\
=&\,
\Big(\sum_{i=0}^3|D_i|\Big) \cdot 0_G
+\Big(\sum_{i=0}^3|D_i|-|G|\Big)G^\ast \nonumber
\end{align}
since ${\mathcal A}$ satisfies (a1), (a2) and $\sum_{i=0}^3|D_i|-|G|+1= 2 \sum_{i=0}^3 |A_i|$. 

 Now, by comparing the two sums  \eqref{Ka:h2} and \eqref{Ka:h3}, it is then easy to see that ${\mathcal D}$ is a symmetric difference family of type $H$ if and only if 
(a4) is satisfied.

{\it Claim 2.} Suppose ${\mathcal A}$ satisfies (a1) and (a2).  Then (a3) holds for ${\mathcal A}$ if and only if ${\mathcal D}$ satisfies (d2). 

By assumption, we have   
\begin{align*}
D_0+D_1-D_2-D_3=\sum_{i=2,3,4,5}A_i-\sum_{i=0,1,6,7}A_i,\\
D_0+D_3-D_1-D_2=\sum_{i=1,2,4,7}A_i-\sum_{i=0,3,5,6}A_i.  
\end{align*}
As ${\mathcal A}$ satisfies (a1) and (a2), we may write 
\[ E_0= \bigcup_{i=2,3,4,5}A_i \mbox{ and } E_1=\bigcup_{i=1,2,4,7}A_i.\]
Our claim is now clear.

\medskip

Now, assume ${\mathcal A}$ is a building family. By Claims 1 and 2, ${\mathcal D}$ is a symmetric difference family of type $H$ that satisfies (d2). Finally, (d1) follows easily from (a2).

Conversely, we assume  ${\mathcal D}=\{D_i:i=0,1,2,3\}$ is a symmetric difference family of type $H$ in $G$ satisfying  (d1) and (d2). By Claims 1 and 2, it suffices to show  that (a1) and (a2) hold. Note that  
\begin{align*}
A_i=\bigcap_{j=0: j\not=i}^3 D_j,\, i=0,1,2,3, 
\end{align*}
and 
\begin{align*}
A_i=\bigcap_{j=0: j\not=i-4}^3 (G^\ast\setminus D_j),\, i=4,5,6,7. 
\end{align*}
It is clear that 
\[
A_i\cap A_j=\bigcap_{h=0}^3 D_h=\emptyset, \quad i,j=0,1,2,3,\, i\not=j, 
\] 
\[
A_i\cap A_j=\bigcap_{h=0}^3 (G^\ast\setminus D_h)=G^\ast \setminus
\bigcup_{h=0}^3 D_h =\emptyset, \quad i,j=4,5,6,7,\, i\not=j, 
\]
and \[
A_i\cap A_j=\Big(G\setminus \bigcup_{h=0: h \not=i}^3 D_h\Big)\cap 
\Big(\bigcap_{h=0: h \not=j}^3 D_h\Big) =\emptyset, \quad i=0,1,2,3,\,  j=4,5,6,7. 
\]
Hence, (a1) holds. Next, for any $x\in G^\ast$, the coefficient of $x$ in $\sum_{i=0}^3D_i$ is either $1$ or $3$ as (d2) is satisfied. If the coefficient of $x$ in $\sum_{i=0}^3D_i$ is $1$, then $x\in \bigcup_{i=4}^7A_i$. 
On the other hand, if the coefficient of $x$ in $\sum_{i=0}^3D_i$ is $3$, then $x\in \bigcup_{i=0}^3A_i$. Hence, $\bigcup_{i=0}^7A_i=G^\ast$ and (a2) is satisfied. \qed



\begin{lemma}\label{lem:key3}
Let ${\mathcal A}=\{ A_i: i=0,\ldots, 7\}$ be  a building family of  an abelian group $G$.  
Then, 
\[ U:=(A_0-A_4)(A_6-A_2)+(A_1-A_5)(A_7-A_3)=-\sum_{i=0}^3A_i.\]
\end{lemma}
\proof  We continue with our notations used in Proposition~\ref{prop:equiv}. Then, we have 
\[
X_0=E_0-\overline{E_0}=\sum_{i=0,1,6,7}A_i-\sum_{i=2,3,4,5}A_i
\mbox{ and }    X_1=E_1-\overline{E_1}=\sum_{i=0,3,5,6}A_i-\sum_{i=1,2,4,7}A_i. 
\]
By Proposition~\ref{prop:equiv}, ${\mathcal D}$ satisfies (d2). Hence, by \eqref{eq:w2}, we have 
\begin{equation}\label{K4:eq}
\sum_{i=0,1}X_i^2=2(|G|-1)\cdot 0_G-2G^\ast. 
\end{equation}
On the other hand, by expanding $X_0^2$ and $X_1^2$, we obtain 
\[\sum_{i=0,1}X_i^2
=2\sum_{i=0}^7A_i^2-2\sum_{i=0}^7A_iA_{i+4}+4U.\]
Then, by (a4) and (a1), we get 
\begin{equation} \label{K4:eq2}
 \sum_{i=0,1}X_i^2=2(|G|-1)\cdot 0_G+2\sum_{i=0}^3A_i-2\sum_{i=4}^7A_i+4U=
 2(|G|-1)\cdot 0_G+4\sum_{i=0}^3A_i-2G^*+4U. 
\end{equation}
Hence, by  \eqref{K4:eq} and \eqref{K4:eq2}, we obtain $U=-\sum_{i=0}^3A_i$. 
\qed

\vspace{0.3cm}

To recap, we have shown that if ${\mathcal D}$ is a symmetric difference family of type $H$ in $G$, then we have 
\begin{itemize}
\item [(a)] ${\mathcal D}$  satisfies (d2) if and only if it satisfies both (d3) and (d4). 
\item [(b)] ${\mathcal D}$  satisfies (d2) if and only if ${\mathcal A}$ defined in Proposition~\ref{prop:equiv} is a building family. 
\end{itemize}
As we will show in the next section, condition (d3) and the notion of building family play a very important role in constructing difference families from a group $G$ to a group $G\times G'$.



\section{Constructing difference families in $G\times G'$ from 
 difference families in $G,G'$}\label{sec:const}

In this section, we exploit a similar idea used in Polhill's paper to construct difference families of type  
$H_4^\ast$ in $G\times G'$ for some suitable group $G'$. Our constructions are based on the existence of difference families in $G$ that satisfy certain conditions. 
First, we consider difference families of type $H$ in $G$ that satisfy (d3).  

\begin{example}\label{ex:conh4ast}
\begin{itemize}
\item[(a)]
Let $v=13$.  Define 
\[
D_0= \{ 1, 2, 3, 5, 6, 9 \},\, 
D_1= \{ 1, 12, 3, 4, 9, 10 \},\,
D_2= \{ 0, 2, 5, 6 \},\,
D_3= \{ 0, 1, 3, 9 \}. 
\]
Then,  ${\mathcal D}=\{D_i:i=0,1,2,3\}$ is a 
difference family of type $H$  in $\Z_{13}$ that satisfies (d3). 
\item[(b)] 
Let $v=37$ and $\omega$ be a suitable primitive element  in $\F_{37}$, define 
\begin{align*}
D_0&\,= \{0\}\cup \{\omega^{i+12j}: i \in \{0,3,7,8,10\},j\in \{0,1,2\}\},\\
D_1&\,= \{\omega^{i+12j}: i \in \{0,1,7,8,10\},j\in \{0,1,2\}\},\\
D_2&\,= \{0\}\cup \{\omega^{i+12j}: i \in \{0,1,5,7,9\},j\in \{0,1,2\}\},\\
D_3&\,=  \{\omega^{i+12j}: i \in \{1,3,4,5,7\},j\in \{0,1,2\}\}. 
\end{align*}
Then, ${\mathcal D}=\{D_i:i=0,1,2,3\}$ is a 
difference family of type $H$  in $\Z_{37}$  that satisfies (d3). 
\end{itemize}
\end{example}

Note that in both examples above, the difference family ${\mathcal D}$ constructed is not symmetric. 

\medskip

For convenience, we will abuse the usage of $\times$ in $\Z[G\times G']$.  For any 
$X=\sum a_g g\in \Z[G]$ and $Y=\sum b_h h\in \Z[G']$, we define $X\times Y=\sum_{g\in G} \sum_{h\in G'} a_gb_h (g,h)\in \Z[G\times G']$. Clearly, for any subsets $X_1, X_2$ of $G$ and  $Y_1, Y_2$ of $G'$, 
$X_1\times Y_1$ and $X_2\times Y_2$ are subsets in $G\times G'$. In $\Z[G\times G']$, 
\[ (X_1\times Y_1) \cdot ( X_2\times Y_2)= (X_1\cdot X_2) \times (Y_1\cdot Y_2).\]

\begin{theorem}\label{le:0022}
Suppose  ${\mathcal D}=\{D_i:i=0,1,2,3\}$ is a difference family of type $H$ in an abelian group $G$   that satisfies (d3). 
If there exists a difference family ${\mathcal S}=\{S_i:i=0,1,2,3\}$ of type $H_4^\ast$ in an abelian group $G$ satisfying (c1), 
 then there exists a  difference family of type $H_4^\ast$  in $G\times \Z_2$. 
\end{theorem}
\proof 
Define 
\begin{align*}
B_0=&\,
 (D_0\times \{0\}) \cup (D_2^c\times \{1\}),\\
B_1=&\,
 (D_1\times \{0\}) \cup (D_3^c\times \{1\}),\\
B_2=&\,
 (S_0\times \{0\}) \cup (S_2\times \{1\}),\\
B_3=&\,
  (S_1\times \{0\}) \cup (S_3\times \{1\}). 
\end{align*}
It is clear that $\sum_{i=0}^3B_iB_i^{(-1)}=X\times 0+Y \times 1$ where
\[ X=  \Big(\sum_{i=0,1}D_iD_i^{(-1)}+\sum_{i=2,3}{D_i^c}{D_i^c}^{(-1)}+\sum_{i=0}^3S_iS_i^{(-1)}\Big) \mbox{ and } \]
\begin{align*}
 Y &\ = 2(|D_0|+|D_1|)G-D_2 D_0^{(-1)}-D_0 D_2^{(-1)}-
D_3 D_1^{(-1)}-D_1 D_3^{(-1)} \\
&\hspace{0.7cm}
+S_0S_2^{(-1)}+S_2S_0^{(-1)}+S_1S_3^{(-1)}+S_3S_1^{(-1)}. 
\end{align*}
As ${\mathcal D}$ is a type $H$ difference family, 
\begin{align*} \sum_{i=0,1}D_iD_i^{(-1)}+\sum_{i=2,3}{D_i^c}{D_i^c}^{(-1)}= & 
\sum_{i=0}^3D_iD_i^{(-1)}+2(|G|-|D_2|-|D_3|)G\\
 = &  \Big(\sum_{i=0}^3|D_i|+ |G|-2|D_2|-2|D_3|\Big)G +|G| \cdot 0_G. \end{align*}
Since ${\mathcal S}$ is of type $H_4^*$ ,   
\[ \sum_{i=0}^3S_iS_i^{(-1)}= (\sum_{i=0}^3|S_i|- |G|-1]G +(|G|+1)\cdot 0_G.\]
Hence, 
\[ X= \left(\sum_{i=0}^3(|D_i|+|S_i|)-2|D_2|-2|D_3|-1\right)G + (2|G|+1)\cdot 0_G.\]
Furthermore, as ${\mathcal D}$ satisfies (d3) and ${\mathcal S}$ satisfies (c1), we have 
\begin{align*}
Y=&\,2(|D_0|+|D_1|)G-\left(\sum_{i=0}^3|D_i|-|G|+1\right)G+
\left(\sum_{i=0}^3|S_i|-|G|\right)G\\
=&\,\left(\sum_{i=0}^3(|D_i|+|S_i|)-1-2|D_2|-2|D_3|\right)G.
\end{align*}
Note that
\[\sum_{i=0}^3 |B_i|-|G\times \Z_2| -1= \sum_{i=0}^{3} (|D|_i+|S_i|)-2|D_2|-2|D_3|-1.\]
Hence, ${\mathcal B}=\{B_i:i=0,1,2,3\}$ forms a difference family of type $H_4^\ast$ in $G\times \Z_2$. 
\qed

\medskip

In order to apply the above result, we need a difference family of type $H_4^\ast$ satisfying  (c1). 
One easy way to construct such families is by using Paley type partial difference sets.

\begin{lemma}\label{lem:Paley4}
Let $P$ be a Paley type partial difference set in  $G$. Then the family
\[ {\mathcal S}=\{  P, P\cup \{0_G\}, P\cup \{0_G\}, G\setminus P\}\]
is a difference family of type $H_4^\ast$ that satisfies (c1).  
\end{lemma}

In particular, there is a difference family of type $H_4^\ast$ satisfying (c1)  for  $|G|=13,37$. At the same time,  we have shown in Example~\ref{ex:conh4ast} that there exist difference families of type $H$ that satisfy (d3). We have thus shown the following:

\begin{proposition} 
There exists a difference family of type $H_4^\ast$  in $\Z_{v}\times \Z_2$ for $v=13,37$. 
\end{proposition}

Next, we consider the case when $G'=\Z_3$.

\begin{theorem}\label{le:00}
Suppose  ${\mathcal D}=\{D_i:i=0,1,2,3\}$ is a symmetric difference family of type $H$ in an abelian group $G$  satisfying (d3).
If there exists a difference family ${\mathcal S}=\{S_i:i=0,1\}$ of type $H_2^\ast$ in an abelian group $G$. 
Then, there exists a  difference family of type $H_4^\ast$  in $G\times \Z_3$. 
\end{theorem}
\proof 
Define 
\begin{align*}
B_0=&\,
 (D_0\times \{1\}) \cup (D_2^c\times \{2\}) 
\cup (S_0\times \{0\}),\\
B_1=&\,
 (D_3\times \{1\}) \cup (D_1^c\times \{2\}) 
\cup (S_1\times \{0\}),\\
B_2=&\,
 (D_0^c\times \{1\}) \cup (D_2\times \{2\}) 
\cup (S_0\times \{0\}),\\
B_3=&\,
  (D_3^c\times \{1\}) \cup (D_1\times \{2\}) 
\cup (S_1\times \{0\}). 
\end{align*}
Since the $D_i$'s are symmetric, it is straightforward to check that 
\begin{align*}
\sum_{i=0}^3B_iB_i^{(-1)}=&
 \Big(\sum_{i=0}^3(D_i^2+{D_i^c}^2)+2\sum_{i=0,1}S_iS_i^{(-1)}\Big)\times \{0\}\\
&+ (D_2^c D_0+D_0^c {D_2}+D_1^c D_3+D_3^c D_1)\times \{1,2\}\\
&+ (S_0^{(-1)}(D_0+D_0^c)+S_1^{(-1)}(D_3+D_3^c)+S_0(D_2+D_2^c)+S_1(D_1+D_1^c))\times \{1\}\\
&+ (S_0(D_0+D_0^c)+S_1(D_3+D_3^c)+S_0^{(-1)}(D_2+D_2^c)+S_1^{(-1)}(D_1+D_1^c))\times \{2\}. 
\end{align*}
By using the assumptions for ${\mathcal D}$ and ${\mathcal E}$, we have  
\[
\sum_{i=0}^3(D_i^2+{D_i^c}^2)+2\sum_{i=0,1}S_iS_i^{(-1)}=
\Big(2\sum_{i=0,1}|S_i|+|G|-1\Big)G^\ast+\Big(4|G|+2\sum_{i=0,1}|S_i|\Big)\cdot 0_G. 
\]
Furthermore, as ${\mathcal D}$ satisfies (d3), we have 
\[
D_2^c D_0+D_0^c {D_2}+D_1^c D_3+D_3^c D_1=\Big(\sum_{i=0}^3|D_i|\Big)G-\Big(\sum_{i=0}^3|D_i|-|G|+1\Big)G=(|G|-1)G. 
\]
Finally, it is clear that 

\pagebreak

\begin{align*}
&\, S_0^{(-1)}(D_0+D_0^c)+S_1^{(-1)}(D_3+D_3^c)+S_0(D_2+D_2^c)+S_1(D_1+D_1^c)\nonumber \\
=&\,S_0(D_0+D_0^c)+S_1(D_3+D_3^c)+S_0^{(-1)}(D_2+D_2^c)+S_1^{(-1)}(D_1+D_1^c)=2\Big(\sum_{i=0,1}|S_i|\Big)G. 
\end{align*}
Summing up, we obtain 
\[ \sum_{i=0}^3B_iB_i^{(-1)}=(4|G|+2\sum_{i=0,1}|S_i|) \cdot 0_{G\times \Z_3}+(2\sum_{i=0,1}|S_i|+|G|-1)(G\times \Z_3)^\ast. \]
This completes the proof of the theorem. 
\qed

\medskip

Unfortunately, condition (d3) is not sufficient to deal with other $G'$. When $G'=\Z_5$, we need a stronger condition on ${\mathcal D}$.  

\begin{theorem}\label{le:002}
Suppose  ${\mathcal D}=\{D_i:i=0,1,2,3\}$ is a symmetric difference family of type $H$ in an abelian group $G$  satisfying  (d2) and (d5). 
Then, there exists a difference family of type $H_4^\ast$  in $G\times \Z_5$. 
\end{theorem}
\proof 
Let $I_0=\{1,2\}$, $I_1=\{3,4\}$, $I_2=\{2,4\}$, and $I_3=\{1,3\}$. Furthermore, we let 
$E_0,E_1$ be as defined in condition (d2). We define 
\begin{align*}
B_0=&\,
 (D_0\times I_0) \cup (D_2^c\times I_1) 
\cup (\overline{E}_0\times \{0\}),\\
B_1=&\,
 (D_3\times I_0) \cup (D_1^c\times I_1) 
\cup (\overline{E}_1\times \{0\}),\\
B_2=&\,
 (D_0^c\times I_2) \cup (D_2\times I_3) 
\cup (E_0\times \{0\}),\\
B_3=&\,
  (D_3^c\times I_2) \cup (D_1\times I_3) 
\cup (E_1\times \{0\}). 
\end{align*}
Note that $D_i$'s and $E_i$'s are  symmetric. For convenience, we write 
\[ \sum_{i=0}^3B_iB_i^{(-1)}=\sum_{i=1}^{6} X_i \mbox{ where } \]
\begin{align}
X_1 =&(D_0^2+D_3^2)\times (I_0I_0^{(-1)}+I_2I_2^{(-1)})+
(D_1^2+D_2^2)\times (I_1I_1^{(-1)}+I_3I_3^{(-1)})\nonumber\\
X_2 = &\,((|D_0|+|D_3|) G-D_0D_2-D_1D_3)\times 
(I_0I_1^{(-1)}+I_1I_0^{(-1)})\nonumber\\
X_3 = &\,((|D_1|+|D_2|) G-D_0D_2-D_1D_3)\times 
(I_2I_3^{(-1)}+I_3I_2^{(-1)})\nonumber\\
X_4 =&\,\Big(2|G|-2\sum_{i=1,2}|D_i|\Big)G\times I_1I_1^{(-1)}+\Big(2|G|-2\sum_{i=0,3}|D_i|\Big)G\times I_2I_2^{(-1)}\nonumber\\
X_5=&\,[\overline{E}_0(D_0+D_2^c)+\overline{E}_1(D_3+D_1^c)+E_0(D_2+D_0^c)+
E_1(D_1+D_3^c)]\times 
\Z_5^\ast\nonumber
\\
X_6 =&\,\sum_{i=0,1}(E_i^2+\overline{E}_i^2)\times 0.\nonumber
\end{align}
Since $I_0I_0^{(-1)}+I_2I_2^{(-1)}=I_1I_1^{(-1)}+I_3I_3^{(-1)}=4\cdot 0+\Z_5^\ast$
and ${\mathcal D}$ forms a difference family of type $H$, 
\[ X_1=\Big(\sum_{i=0}^3 D_i^2\Big) \times (4\cdot 0 +\Z_5^\ast)
=\,[\Big(\sum_{i=0}^3|D_i|\Big)\cdot 0_G+\Big(\sum_{i=0}^3|D_i|-|G|\Big)G^\ast]\times (4 \cdot 0 +\Z_5^\ast). \]
Note that  $|D_0|+|D_3|=|D_1|+|D_2|$, and $ I_0I_1^{(-1)}+I_1I_0^{(-1)}+I_2I_3^{(-1)}+I_3I_2^{(-1)}=4\Z_5^\ast$. Moreover, ${\mathcal D}$ satisfies (d3)  by Proposition~\ref{prop:key1}. Hence,
\[ 2(D_0D_2+D_1D_3)= (\sum_{i=0}^3|D_i|-|G|+1)G.\]
 Therefore, 
\[ X_2+X_3= \frac{1}{2}(|G|-1)G\times 
(I_2I_3^{(-1)}+I_3I_2^{(-1)})=2(|G|-1)(G\times \Z_5^\ast). \]
Again, as $|D_0|+|D_3|=|D_1|+|D_2|$, 
\[ 2|G|-2\sum_{i=1,2}|D_i|= 2|G|-2\sum_{i=0,3}|D_i|=2|G|-\sum_{i=0}^3 |D_i|.\]
Therefore, 
\[ X_4= \Big(2|G|-\sum_{i=0}^3|D_i|\Big)G\times (I_1I_1^{(-1)}+I_2I_2^{(-1)})= \Big(2|G|-\sum_{i=0}^3|D_i|\Big)G\times (4\cdot 0+\Z_5^\ast). \]
For $X_5$, since 
 ${\mathcal D}$ satisfies (d2), it follows from Lemma~\ref{lem:key2} that 
\[ 
\overline{E}_0(D_0+D_2^c)+\overline{E}_1(D_3+D_1^c)+E_0(D_2+D_0^c)+
E_1(D_1+D_3^c)=T_0=(|G|-1)\cdot 0_G+(2|G|-1)G^\ast .
\]
Therefore, 
\[ X_5=[ (|G|-1)\cdot 0_G+(2|G|-1)G^\ast] \times 
\Z_5^\ast.\]
Finally, as $E_0, E_1, \overline{E}_0, \overline{E}_2$ forms a difference family of type $H_4^*$, we get 
\[ X_6=2(|G|-1)\cdot 0_G+(|G|-3)(G^\ast \times 
\Z_5^\ast).\]
Summing up,  we obtain 
$\sum_{i=0}^3B_iB_i^{(-1)}=(10|G|-2)\cdot 0_{G\times \Z_5}+(5|G|-3)(G\times \Z_5)^\ast$.  
This completes the proof of the theorem. 
\qed

\medskip

It seems difficult to generalize the above constructions for other $G'$. However, instead of constructing symmetric difference families of type $H_4^*$,  it is possible to construct symmetric difference families of type $H_8^\ast$ under stronger assumptions. In the previous examples, we simply make use of the sets $D_i$ or $D_i^c$. But in the next construction, we need to make use of a building family in $G'$.
 
\begin{theorem}\label{le:01}
Suppose ${\mathcal D}=
\{D_i:i=0,1,2,3\}$ is a symmetric difference family in an abelian group $G$  satisfying the (d2). 
If there exists a building family ${\mathcal B}=\{B_i:i=0,1,\ldots,7\}$ in an abelian group $G'$, then, there exists a symmetric difference family of type $H_8^\ast$ in $G\times G'$. 
\end{theorem}
{\bf Proof.\, }  As  ${\mathcal D}$ satisfies (d2), we let $E_0,E_1$ be as defined in (d2). We define 
\begin{align*}
C_0=&\,
 (D_0\times B_0) \cup (D_1\times B_1) \cup 
(D_2^c\times B_2) \cup (D_3^c\times B_3) \\
&\, \cup 
(D_0^c\times B_4) \cup (D_1^c\times B_5) \cup 
(D_2\times B_6) \cup (D_3\times B_7) \cup (\overline{E_0}\times \{0_{G'}\}),\\
C_1=&\,
 (D_0\times B_1) \cup (D_1\times B_4) \cup 
(D_2^c\times B_3) \cup (D_3^c\times B_6) \\
&\, \cup 
(D_0^c\times B_5) \cup (D_1^c\times B_0) \cup 
(D_2\times B_7) \cup (D_3\times B_2) \cup (\overline{E_1}\times \{0_{G'}\}),\\
C_2=&\,
 (D_0\times B_2) \cup (D_1\times B_3) \cup 
(D_2^c\times B_0) \cup (D_3^c\times B_1) \\
&\, \cup 
(D_0^c\times B_6) \cup (D_1^c\times B_7) \cup 
(D_2\times B_4) \cup (D_3\times B_5) \cup (\overline{E_0}\times \{0_{G'}\}),\\
C_3=&\,
 (D_0\times B_3) \cup (D_1\times B_6) \cup 
(D_2^c\times B_1) \cup (D_3^c\times B_4) \\
&\, \cup 
(D_0^c\times B_7) \cup (D_1^c\times B_2) \cup 
(D_2\times B_5) \cup (D_3\times B_0) \cup (\overline{E_1}\times \{0_{G'}\}),\\
C_4=&\,
 (D_0\times B_4) \cup (D_1\times B_5) \cup 
(D_2^c\times B_6) \cup (D_3^c\times B_7) \\
&\, \cup 
(D_0^c\times B_0) \cup (D_1^c\times B_1) \cup 
(D_2\times B_2) \cup (D_3\times B_3) \cup (E_0\times \{0_{G'}\}),\\
C_5=&\,
 (D_0\times B_5) \cup (D_1\times B_0) \cup 
(D_2^c\times B_7) \cup (D_3^c\times B_2) \\
&\, \cup 
(D_0^c\times B_1) \cup (D_1^c\times B_4) \cup 
(D_2\times B_3) \cup (D_3\times B_6) \cup (E_1\times \{0_{G'}\}),\\
C_6=&\,
 (D_0\times B_6) \cup (D_1\times B_7) \cup 
(D_2^c\times B_4) \cup (D_3^c\times B_5) \\
&\, \cup
(D_0^c\times B_2) \cup (D_1^c\times B_3) \cup 
(D_2\times B_0) \cup (D_3\times B_1) \cup (E_0\times \{0_{G'}\}),\\
C_7=&\,
 (D_0\times B_7) \cup (D_1\times B_2) \cup 
(D_2^c\times B_5) \cup (D_3^c\times B_0) \\
&\, \cup
(D_0^c\times B_3) \cup (D_1^c\times B_6) \cup 
(D_2\times B_1) \cup (D_3\times B_4) \cup (E_1\times \{0_{G'}\}). 
\end{align*}
It is clear that each $C_i$ is symmetric. 
We now prove that 
\begin{equation}\label{eq:Bilast}
\sum_{i=0}^7C_i^2=4(|G|\cdot |G'|-1) \cdot 0_{G\times G'}+2(|G|\cdot |G'|-3) (G\times G')^\ast. 
\end{equation}
To simplify notations, we set 
\begin{align*}
P=&\, \sum_{i=0}^7B_i^2-\sum_{i=0}^7 B_i B_{i+4},\\
Q=&\, \Big(\sum_{i=2,3,4,5}B_i\Big)\Big(\sum_{i=0,1,6,7}B_i\Big)+
\Big(\sum_{i=0,3,5,6}B_i\Big)\Big(\sum_{i=1,2,4,7}B_i\Big),\\
R=&\, \Big(\sum_{i=2,3,4,5}B_i\Big)^2+\Big(\sum_{i=0,1,6,7}B_i\Big)^2
+\Big(\sum_{i=0,3,5,6}B_i\Big)^2+\Big(\sum_{i=1,2,4,7}B_i\Big)^2. 
\end{align*} 
As before, we let $W,T_0, T_1, T_2, T_3$ be as defined in Section~\ref{sec:proper} and 
\[ U:=(B_0-B_4)(B_6-B_2)+(B_1-B_5)(B_7-B_3). \]

 Then, we have 
\begin{align*}
\sum_{i=0}^7B_i^2=&\, 2\Big(\sum_{i=0}^3D_i^2\Big)\times P+8W\times U+2 \Big(\sum_{i=0}^3|D_i|\Big)G\times Q\\
&\, +\big(2|G|-\sum_{i=0}^3|D_i|\Big) G\times R+2\Big(\sum_{i=0,1}E_i^2+\sum_{i=0,1}\overline{E}_i^2\Big)\times 0_{G'}\\
&\, +2[T_0\times (B_0+B_2)+T_1\times (B_1+B_3)+
T_2\times (B_4+B_6)+T_3\times (B_5+B_7)]. 
\end{align*}
By 
(a4), we have 
\[
2\Big(\sum_{i=0}^3D_i^2\Big)\times P
=[2\Big(\sum_{i=0}^3|D_i|\Big) \cdot 0_G+2\Big(\sum_{i=0}^3|D_i|-|G|\Big)G^\ast]\times 
[(|G'|-1)\cdot 0_{G'}+\sum_{i=0}^3B_i-\sum_{i=4}^7B_i]. 
\]
Next, 
by Proposition~\ref{prop:key1} and Lemma~\ref{lem:key3}, we have 
\[
8W\times U=-4\Big(\sum_{i=0}^3|D_i|-|G|+1\Big)G\times \Big(\sum_{i=0}^3B_i\Big). 
\]
Furthermore, by (a3), we have 
\begin{align*}
&\,2 \Big(\sum_{i=0}^3|D_i|\Big)G\times Q+\Big(2v-\sum_{i=0}^3|D_i|\Big) G\times R\\
=&\,2\Big(\sum_{i=0}^3|D_i|+|G|(|G'|-3)\Big) G\times (G')^\ast
+2(|G'|-1)\Big(2|G|-\sum_{i=0}^3|D_i|\Big)(G\times 0_{G'}). 
\end{align*}
Moreover, by (d2), 
\[
2[\sum_{i=0,1}E_i^2+\sum_{i=0,1}\overline{E}_i^2]\times 0_{G'}=[4(|G|-1)\cdot 0_G+2(|G|-3)G^\ast]\times 0_{G'}. 
\]
Finally, 
by Lemma~\ref{lem:key2}, we have  
\begin{align*}
&\, 2[T_0\times (B_0+B_2)+T_1\times (B_1+B_3)+
T_2\times (B_4+B_6)+T_3\times (B_5+B_7)]\\
=&\, [(2(|G|-1)\cdot 0_G+(2|G|-1)G^\ast]\times \Big(\sum_{i=0}^3B_i\Big)+2[3(|G|-1)\cdot 0_G+(2|G|-3)G^\ast]\times \Big(\sum_{i=4}^7B_i\Big). 
\end{align*}
Summing up, we obtain \eqref{eq:Bilast}. 
This completes the proof of the theorem. 
\qed 

\medskip

In view of Theorems~\ref{le:0022}, \ref{le:00}, \ref{le:002} and  \ref{le:01}, we have then completed the proof of Theorem~\ref{thm:consttypeH2}. 

\section{Symmetric difference families of type $H$ that satisfy (d1) and (d2)
}\label{sec:existsDF}

To apply the constructions we obtain in the previous section, it remains to find groups that contain symmetric difference families of type $H$ satisfying (d1) and (d2). 
In \cite{XL91}, the authors constructed difference families of type $H$ in $(\F_{q^2},+)$ for $q\equiv 1\,(\mod{4})$. Our first result is to show that those difference families actually satisfy (d1) and (d2). Let us recall the construction.


Let $q=4m+1$ be a prime power and let $\omega$ be a primitive 
element of $\F_{q^2}$. Let $N$ be a divisor of $q^2-1$. We define 
\[
C_i^{(N,q^2)}=\omega^i\langle \omega^{N}\rangle, \, \, \, \, i=0,1,\ldots,N-1,
\] 
and 
\[
D_0:=\Big(\bigcup_{\ell=0}^{m-1}\bigcup_{u=1}^3 C_{4\ell+(2m+1)u}^{(8m+4,q^2)}\Big)
\cup \Big(\bigcup_{j=0}^m C_{4j-2}^{(8m+4,q^2)}\Big), \, \, 
D_i=\omega^{(2m+1)i}D_0,\, \,  i=1,2,3. 
\]

It has been shown in \cite{XL91} that ${\mathcal D}=\{D_i:i=0,1,2,3\}$ forms a
difference family of type $H$ in $(\F_{q^2},+)$. 
It is clear that each $D_i$ is symmetric and $0\not\in \bigcup_{i=0}D_i$. Therefore, (d1) is satisfied. 
Furthermore, we have 
\begin{align*}
D_0+D_1-D_2-D_3=&\, \Big(\sum_{\ell=0}^{m-1}\sum_{u=2,3}C_{4\ell+(2m+1)u}^{(8m+4,q^2)}+
\sum_{\ell=0}^{m}\sum_{u=0,1}C_{4j-2+(2m+1)u}^{(8m+4,q^2)}
\Big)\\
&\hspace{1cm}-\Big(\sum_{\ell=0}^{m-1}\sum_{u=0,1}C_{4\ell+(2m+1)u}^{(8m+4,q^2)}+
\sum_{j=0}^{m}\sum_{u=2,3}C_{4j-2+(2m+1)u}^{(8m+4,q^2)}
\Big) \\
=&\, (C_2^{(4,q^2)}+C_{2m+3}^{(4,q^2)})-(C_0^{(4,q^2)}+C_{2m+1}^{(4,q^2)}), \\
D_0+D_3-D_1-D_2=&\, \Big(\sum_{\ell=0}^{m-1}\sum_{u=1,2}
C_{4\ell+(2m+1)u}^{(8m+4,q^2)}+
\sum_{j=0}^{m}\sum_{u=0,3}C_{4j-2+(2m+1)u}^{(8m+4,q^2)}
\Big)\\
&\hspace{1cm}-\Big(\sum_{\ell=0}^{m-1}\sum_{u=0,3}C_{4\ell+(2m+1)u}^{(8m+4,q^2)}+
\sum_{j=0}^{m}\sum_{u=1,2}C_{4j-2+(2m+1)u}^{(8m+4,q^2)}
\Big) \\
=&\, (C_2^{(4,q^2)}+C_{2m+1}^{(4,q^2)})-(C_0^{(4,q^2)}+C_{2m+3}^{(4,q^2)}). 
\end{align*}
It is known (cf.~\cite{Y92}) that the set of $E_0=C_2^{(4,q^2)}\cup C_{2m+3}^{(4,q^2)}$, 
$\overline{E}_0=C_0^{(4,q^2)}\cup C_{2m+1}^{(4,q^2)}$, 
$E_1=C_2^{(4,q^2)}\cup C_{2m+1}^{(4,q^2)}$, 
$\overline{E}_1=C_0^{(4,q^2)}\cup C_{2m+3}^{(4,q^2)}$ forms a difference family 
of type $H_4^\ast$ in $(\F_{q^2},+)$. Hence, ${\mathcal D}$ satisfies (d2). 
Note that all $|D_i|$'s are equal. Therefore, (d5) is satisfied. As ${\mathcal D}$ is also symmetric, it follows from Proposition 2.3 that (d3) and (d4) are also satisfied. 
We have thus proved the following: 

\begin{theorem}\label{thm:ex1}
Let $q\equiv 1\,(\mod{4})$ be a prime power. Then, 
there exists a symmetric difference family of type $H$ in $(\F_{q^2},+)$ that satisfies (d1)-(d5).
\end{theorem}



Other constructions of symmetric  difference family  ${\mathcal D}$ of type $H$ in $G=(\F_q^4,+)$  
are found in \cite{C97,X92,WX97,XC96}.

\begin{proposition}\label{prop:exists}
Let $G=(\F_9,+)$ or $G=(\F_q^4,+)$ where $q$ is an odd prime power. There exists a symmetric difference family  ${\mathcal D}$ of type $H$ in $G$ such that the following conditions are satisfied. 
\begin{itemize}
\item[(i)] For any $i=0,1,2,3$, $|D_i|=\frac{|G|-\sqrt{|G|}}{2}$. 
\item [(ii)]  For any nontrivial character $\psi$ of $G$, exactly one of the values $\psi(D_i)$, $i=0,1,2,3$, 
is nonzero; and is equal to $\pm \sqrt{|G|}$. 
\item [(iii)] $H_0+H_1+H_4+H_5=G+1\cdot 0_G$ and $0_G\not \in \bigcup_{i=0}^3D_i$ or 
$0_G\in \bigcap_{i=0}^3D_i$ if we set 
\begin{align*}
&H_0=D_0\cap D_1,H_1=D_0^c\cap D_1^c,H_2=D_0\cap D_1^c,H_3=D_0^c\cap D_1,\\
&H_4=D_2\cap D_3,H_5=D_2^c\cap D_3^c,H_6=D_2\cap D_3^c,H_7=D_2^c\cap D_3. 
\end{align*}
\end{itemize}
Furthermore, we may assume $0_G\not \in \bigcup_{i=0}^3D_i$ if $G=(\F_q^4,+)$ and $0_G\in \bigcap_{i=0}^3D_i$ if $G=(\F_9,+)$. 
 \end{proposition}

\begin{remark}\label{rem:fin}
\begin{itemize}
\item [(a)] Note that by using  conditions (i) and (ii) above and   the Fourier transform, it is clear that ${\mathcal D}$ satisfies (d3) as we have
\[ 2(D_0D_2+D_1D_3)=(\sqrt{|G|}-1)^2G=(\sum_{i=0}^3|D_i|-|G|+1)G.  \]
\item [(b)] ${\mathcal D}$ satisfies (d4) if condition (iii) is satisfied since $\sum_{i=0}^3D_i=2(H_0+H_4)+H_2+H_3+H_6+H_7=2(H_0+H_4)+G^\ast$. 
\item [(c)] Therefore, by Proposition~\ref{prop:key1},  ${\mathcal D}$ satisfies (d2) if ${\mathcal D}$ satisfies conditions (i), (ii)  and (iii) above.
\end{itemize}
\end{remark} 

In view of the above remark, we obtain another symmetric difference families of type $H$ that satisfy (d2). 
However, we do not obtain new parameters from these families. To construct one with new parameters, we employ a standard technique.


Suppose ${\mathcal D}=\{D_i:i=0,1,2,3\}$ and ${\mathcal F}=\{F_i:=i=0,1,2,3\}$ are symmetric difference families of type $H$ in abelian groups $G$ and $N$ respectively. 
If both ${\mathcal D}$ and ${\mathcal F}$
satisfy conditions (i), (ii) and (iii) in Proposition~\ref{prop:exists}, we then define the following:
\begin{align*}
B_0&\,=H_1\times F_0+H_0\times F_0^c+H_3\times F_2+H_2\times F_2^c\\
B_1&\,=H_1\times F_1+H_0\times F_1^c+H_2\times F_3+H_3\times F_3^c\\
B_2&\,=H_6\times F_0+H_7\times F_0^c+H_5\times F_2+H_4\times F_2^c\\
B_3&\,=H_7\times F_1+H_6\times F_1^c+H_5\times F_3+H_4\times F_3^c.
\end{align*}
Clearly, ${\mathcal B}$ satisfies condition (i). As shown in (cf. \cite{P10}), ${\mathcal B}=\{B_i:i=0,1,2,3\}$ is indeed a symmetric difference family of type $H$ in $G\times N$ satisfying condition (ii). We now check that  ${\mathcal B}$ satisfies condition (iii) under the assumption that $0_N\not \in \bigcup_{i=0}^3F_i$: 
\begin{align}
&\,(B_0\cap B_1)+(B_0^c\cap B_1^c)+ (B_2\cap B_3)+ (B_2^c\cap B_3^c)\nonumber\\
=&\,H_1\times (F_0\cap F_1)+H_0\times (F_0^c\cap F_1^c)+ 
H_3\times (F_2\cap F_3^c)+ H_2\times (F_3\cap F_2^c)\nonumber\\
&\,+H_1\times (F_0^c\cap F_1^c)+H_0\times (F_0\cap F_1)+ 
H_3\times (F_2^c\cap F_3)+ H_2\times (F_3^c\cap F_2)\nonumber\\
&\,+H_6\times (F_0\cap F_1^c)+H_7\times (F_0^c\cap F_1)+ 
H_5\times (F_2\cap F_3)+ H_4\times (F_3^c\cap F_2^c)\nonumber\\
&\,+H_6\times (F_0^c\cap F_1)+H_7\times (F_0\cap F_1^c)+ 
H_5\times (F_2^c\cap F_3^c)+ H_4\times (F_3\cap F_2)\nonumber\\
=&\,(H_0+H_1)\times ((F_0\cap F_1)\cup (F_0^c\cap F_1^c))+ 
(H_6+H_7)\times ((F_0\cap F_1^c)+  (F_0^c\cap F_1))\nonumber\\
&\,+(H_2+H_3)\times ((F_2^c\cap F_3)\cup (F_3^c\cap F_2))+ 
(H_4+H_5)\times ((F_2\cap F_3)+ (F_3^c\cap F_2^c)).\label{eq:sumeqG} 
\end{align}
Since $H_0+H_1+H_4+H_5=G+1\cdot 0_G$ and $\sum_{i=0}^3H_i=\sum_{i=4}^7H_i=G$, 
we have 
\[ H_6+H_7=H_0+H_1-1\cdot 0_G \mbox{ and }H_2+H_3=H_4+H_5-1\cdot 0_G. \]
It then follows from \eqref{eq:sumeqG} that we have 
\begin{align*}
&\,(B_0\cap B_1)\cup (B_0^c\cap B_1^c)\cup (B_2\cap B_3)\cup (B_2^c\cap B_3^c)\\
&\,=(H_0+H_1+H_4+H_5)\times N- 0_G\times 
((F_2^c\cap F_3)\cup (F_3^c\cap F_2)\cup (F_0^c\cap F_1)\cup (F_1^c\cap F_0)) \\
&\, =(G+1\cdot 0_G)\times N - 0_G\times N^*=G\times N+ 1\cdot 0_{G\times N}
\end{align*}
as $(F_2^c\cap F_3)\cup (F_3^c\cap F_2)\cup (F_0^c\cap F_1)\cup (F_1^c\cap F_0)=N^\ast$. Finally, we need to compute $\sum_{i=0}^3B_i$. Observe that the coefficient of $0_{G\times N}$ in $\sum_{i=0}^3B_i$ is zero if $0_G\not \in \bigcup_{i=0}^3D_i$ and four if $0_G\in \bigcap_{i=0}^3D_i$. That means
$0_G\notin \bigcup_{i=0}^3B_i$ or $0_G\in \bigcap_{i=0}^3B_i$. Therefore, ${\mathcal B}$ satisfies (iii). 

By Remark~\ref{rem:fin}, the symmetric difference family ${\mathcal B}$ satisfies (d2). However, it is not clear if 
it satisfies (d1). But in view of Proposition~\ref{prop:comp}, we conclude that either ${\mathcal B}$ or 
${\mathcal B'}=\{ G\times N - B_i: i=0,1,2,3\}$ is a symmetric difference family of type $H$ that satisfies (d1) and (d2). Finally, as ${\mathcal B}$ and ${\mathcal B'}$ satisfy (i), (d5) holds. To summarize, we have shown the following:


\begin{proposition}\label{prop:recur}
Suppose there exist symmetric difference families ${\mathcal D}=\{D_i:i=0,1,2,3\}$ and ${\mathcal F}=\{F_i:i=0,1,2,3\}$ of type $H$ satisfying conditions (i), (ii) and (iii) in abelian groups $G$ and $N$, respectively. Furthermore, we assume $0_N\not \in \bigcup_{i=0}^3F_i$. Then, there exists a 
symmetric difference family ${\mathcal B}=\{B_i:i=0,1,2,3\}$ of type $H$ satisfying conditions (i), (ii) and (iii) in $G\times N$. Moreover,  either ${\mathcal B}$ or ${\mathcal B}'=\{G\setminus B_i:i =0,1,2,3\}$ is a symmetric difference family of type $H$ that satisfies (d1)-(d5).
\end{proposition} 

By applying Propositions ~\ref{prop:exists} and ~\ref{prop:recur}  repeatedly, we obtain the following theorem.  

\begin{theorem}\label{thm:ex2}
Let $p_1,p_2,\ldots, p_s$ be distinct primes and 
$G=N\times \Z_{p_1}^{4t_1}\times \Z_{p_2}^{4t_2}\times \cdots  \Z_{p_s}^{4t_s}$, where $N=\{0\}$ or 
$\Z_3^2$. Then, there exists a symmetric difference family of type $H$ 
in $G$ satisfying (d1)-(d5). 
\end{theorem}

\section{Conclusion}
In this section, we prove Theorem~\ref{thm:existHada}.  Recall that 
\begin{align*}
\Phi_1=&\,\{q^2:q\equiv 1\,(\mod{4})\mbox{ is a prime power}\},\\
\Phi_2=&\{n^4\in \N:n\equiv 1\,(\mod{2})\} \cup \{9n^4\in \N:n\equiv 1\,(\mod{2})\}, \\
 \Phi_3=&\{5\},
\Phi_4=\{13,37\}.
\end{align*}

Note that for $v\in \Phi_1\cup \Phi_2$,  it follows from Theorems~\ref{thm:ex1} and  \ref{thm:ex2} that there exists a symmetric difference family of type $H$ that satisfies (d1)-(d5). 
For $v\in \Phi_3\cup \Phi_4$, difference families constructed explicitly in Examples~\ref{exam:v5} and \ref{ex:conh4ast} satisfy (d3). 

On the other hand, it is known that a Paley type partial difference set $P$ exists in an abelian group of order $v$ for any $v\in  \Phi_1 \cup \Phi_2 \cup \Phi_3 \cup \Phi_4$ (see, e.g., \cite{P10}). By Lemma~\ref{lem:Paley4}, there exists a difference family of type $H_4^\ast$ satisfying (c1) in an abelian group of order $v$. Also, $\{P, \bar{P}\}$ is a difference family of type $H_2^*$. Therefore, by applying Theorems~\ref{le:0022}, \ref{le:00} and \ref{le:002},  we obtain the following:


\begin{theorem}\label{thm:main001}
There exists a difference family of type $H_4^\ast$ 
in an abelian group of order $2v$ for any $v\in \Phi_1 \cup \Phi_2 \cup \Phi_3 \cup \Phi_4$. 
\end{theorem}

\begin{theorem}\label{thm:main01}
There exists a difference family of type $H_4^\ast$ 
in an abelian group of order $3v$ for any $v\in \Phi_1 \cup \Phi_2 \cup \Phi_3$. 
\end{theorem}

\begin{theorem}\label{thm:main02}
There exists a  difference family of type $H_4^\ast$ 
in an abelian group of order $5v$ for any $v\in \Phi_1 \cup \Phi_2 \cup \Phi_3$. 
\end{theorem}

Note that in Theorems~\ref{thm:main01} and \ref{thm:main02}, $v\notin \Phi_4$ as the family constructed in Example~\ref{ex:conh4ast} is not symmetric.
The next result is now obvious in view of Theorem~\ref{le:01}.

\begin{theorem}\label{thm:main03}
There exists a symmetric difference family of type $H_8^\ast$ 
in an abelian group of order $uv$ for any $u\in \Phi_1 \cup \Phi_2$ and $v\in \Phi_1 \cup \Phi_2 \cup \Phi_3$. 
\end{theorem}

Finally, by using the Wallis-Whiteman array together with the difference families in 
Theorems~\ref{thm:main001}, \ref{thm:main01} and ~\ref{thm:main02}, we obtain (1), (2) 
and (3) in Theorem~\ref{thm:existHada}. By using the Kharaghani array  together with the difference families in 
Theorem~\ref{thm:main03}, we obtain (4) in Theorem~\ref{thm:existHada}. 

\end{document}